\documentclass[a4paper, 10pt]{article}
\usepackage{latexsym,amssymb}
\usepackage[english]{babel}
\usepackage{bm}
\usepackage{amsfonts, amsthm}
\usepackage{amsmath}
\usepackage{mathrsfs}
\usepackage{multirow}
\usepackage{bm}
\usepackage{pgf,tikz}
\usepackage{rotating}
\usepackage{calculator}

\usepackage{xcolor}

\newcommand{\comment}[1]{}

\newcommand{\Keywords}[1]{\par\noindent{\bfseries Keywords}: #1}
\newcommand{\MSC}[1]{\par\noindent{\bfseries MSC}: #1}

\newcommand{\PaperTitle}[1]{#1}

\newcommand{\JournalName}[1]{\emph{#1}}

\def\Z{\mathbb{Z}}
\def\N{\mathbb{N}}

\usepackage{amsthm}
\theoremstyle{plain}

\newtheorem*{Theorem*}{Theorem 1}

\newtheorem*{Lemma*}{Lemma}
\newtheorem{Definition}{Definition}[section]
\newtheorem{Remark}[Definition]{Remark}
\newtheorem{Theorem}[Definition]{Theorem}
\newtheorem{Proposition}[Definition]{Proposition}
\newtheorem{Example}[Definition]{Example}
\newtheorem{Lemma}[Definition]{Lemma}
\newtheorem{Problem}[Definition]{Problem}
\newtheorem{Corollary}[Definition]{Corollary}

\title{
Factorizing the Rado graph and infinite complete graphs
}
\author{Simone Costa\thanks{DII/DICATAM - Sez. Matematica, Universit\`a degli Studi di Brescia, Via
Branze 38, I-25123 Brescia, Italy. email: simone.costa@unibs.it}\ \
and Tommaso Traetta%
\thanks{
DICATAM - Sez. Matematica, Universit\`a degli Studi di Brescia,
Via Valotti 9, I-25123 Brescia, Italy. email: tommaso.traetta@unibs.it}}
\begin{document}
\maketitle

\begin{abstract}
Let $\mathcal{F}=\{F_{\alpha}: \alpha\in \mathcal{A}\}$ be a family of infinite graphs, together with $\Lambda$. The Factorization Problem $FP(\mathcal{F}, \Lambda)$
asks whether $\mathcal{F}$ can be realized as a factorization of $\Lambda$, namely, whether
there is a factorization $\mathcal{G}=\{\Gamma_{\alpha}: \alpha\in \mathcal{A}\}$ of 
$\Lambda$ such that each $\Gamma_{\alpha}$ is a copy of $F_{\alpha}$.

We study this problem when $\Lambda$ is either the Rado graph $R$ 
or the complete graph $K_\aleph$ of infinite order $\aleph$. 
When $\mathcal{F}$ is a countable family,
we show that $FP(\mathcal{F}, R)$ is solvable if and only if 
each graph in $\mathcal{F}$ has no finite dominating set.
%
%
%
We also prove that $FP(\mathcal{F}, K_\aleph)$ admits a solution 
whenever the cardinality $\mathcal{F}$ coincide with the order and the domination numbers of its graphs. 

For countable complete graphs, we show some non existence results when
the domination numbers of the graphs in $\mathcal{F}$ are finite. More precisely,
we show that there is no factorization of $K_\N$ into copies of a $k$-star (that is, 
the vertex disjoint union of $k$ countable stars)  when $k=1,2$, whereas it exists when
$k\geq 4$, leaving the problem open for $k=3$.

Finally, we determine sufficient conditions for 
the graphs of a decomposition to be arranged into
resolution classes.


\end{abstract}
\Keywords{Factorization Problem, Resolution Problem, Rado Graph, Infinite Graphs.}
\MSC{05C63, 05C70}
\section{Introduction}
We assume that the reader is familiar
with the basic concepts in (infinite) graph theory, and refer to 
\cite{Die17} for further details.

In this paper all graphs will be simple, namely, without multiple edges or loops. 
As usual, we denote by $V(\Lambda)$ and $E(\Lambda)$ the vertex set and the edge set of a simple graph $\Lambda$, respectively. We say that $\Lambda$ is finite (resp. infinite) 
if its vertex set is so, and refer to the cardinality of $V(\Lambda)$ and $E(\Lambda)$ 
as the order and the size of $\Lambda$, respectively.
Note that in the finite case $|E(\Lambda)|\leq \binom{|V(\Lambda)|}{2}$, whereas if $\Lambda$ is infinite, then its order, which is a cardinal number, is greater than or equal to its size.
We use the notation $K_v$ for any complete graph of order $v$, and denote by $K_V$ the complete graph whose vertex set is $V$.

Given a subgraph $\Gamma$ of a simple graph $\Lambda$, we denote by
$\Lambda\setminus\Gamma$ the graph obtained from $\Lambda$ by deleting
the edges of $\Gamma$.  
If $\Gamma$ contains all possible edges of $\Lambda$ joining any two of its vertices,
then $\Gamma$ is called an induced subgraph of $\Lambda$ 
(in other words, an induced subgraph is obtained by vertex deletions only).
Instead, if $V(\Gamma) = V(\Lambda)$, then $\Gamma$ is
called a spanning subgraph or a factor of $\Lambda$
(hence, a factor is obtained by edge deletions only). 
If $\Gamma$ is also $h$-regular, then we speak of an $h$-factor.
We recall that a set $D$ of vertices of $\Lambda$ is dominating if all other vertices of $\Lambda$ are adjacent to some vertex of $D$. The minimum size of a dominating set of 
$\Lambda$ is called the domination number of $\Lambda$. Finally, we say that 
$\Lambda$ is locally finite if its vertex degrees are all finite.

A decomposition of $\Lambda$ is a set $\mathcal{G}=\{\Gamma_1,\ldots, \Gamma_n\}$ of subgraphs  
of $\Lambda$ whose edge-sets partition $E(\Lambda)$. 
If the graphs $\Gamma_i$ are all isomorphic to a given subgraph $\Gamma$ of $\Lambda$, then we speak of a 
$\Gamma$-decomposition of $\Lambda$. When $\Gamma$ and $\Lambda$ are both complete graphs, we obtain $2$-designs. More precisely, a $K_k$-decomposition of $K_v$ is equivalent to a $2$-$(v,k,1)$ design. 

Classically, the graphs $\Gamma_i$ and $\Lambda$ are taken to be finite, 
and the same usually holds for the parameters $v$ and $k$ of a $2$-designs. 
However, there has been considerable interest in designs on a infinite set of $v$ points, 
mainly when $k=3$. In this case,  we obtain infinite Steiner triple systems whose first 
explicit constructions were given in \cite{GrGrPh87, GrGrPh91}.
Further results concerning the existence of rigid, sparse, and perfect countably Steiner triple systems can be found in \cite{CaWe12, ChiGrGrWe10, Fr94}. Results showing that any Steiner system can be extended are given in \cite{BeuCa94, Qua80}. The existence of large sets of
Steiner triple systems for every infinite $v$ (and more generally, of infinite Steiner systems)
can be found in \cite{Ca95}. Also, 
infinite versions of topics in finite geometry, including infinite Steiner triple systems and infinite perfect codes are considered in \cite{Ca84}.
A more comprehensive list of results on infinite designs can be found in \cite{DaHoWe14}.

When each graph of a decomposition $\mathcal{G}$ of $\Lambda$ is a factor (resp. $h$-factor), 
we speak of a factorization (resp. $h$-factorization) of $\Lambda$. Also, when the factors of
$\mathcal{G}$ are all isomorphic to $\Gamma$, we speak of a 
$\Gamma$-factorization of $\Lambda$.
A factorization of $K_v$ into factors whose components are copies of $K_k$ is equivalent to
a resolvable $2$-$(v,k,1)$ design.

In this paper, we consider the Factorization Problem for infinite graphs, which is here stated in its most general version

\begin{Problem}
\label{factorization problem}
Let $\Lambda$ be a graph of order $\aleph$ and let $\mathcal{F}=\{F_{\alpha}: \alpha\in \mathcal{A}\}$ be a family of (non-empty) infinite graphs, 
not necessarily distinct, each of which has order $\aleph$, with $\aleph \geq |\mathcal{A}|$.

The Factorization Problem $FP(\mathcal{F}, \Lambda)$ asks for a factorization 
$\mathcal{G}=\{\Gamma_{\alpha}:\alpha\in \mathcal{A}\}$ of $\Lambda$ such that $\Gamma_{\alpha}$ is isomorphic to $F_{\alpha}$, 
for every $\alpha\in\mathcal{A}$. 
If $\Lambda$ is the complete graph of order $\aleph$, we simply write $FP(\mathcal{F})$. If in addition to this each $F_\alpha$ is isomorphic to a given graph $F$ and $|\mathcal{A}|=\aleph$, we write $FP(F)$.
\footnote{Since in this case the factorization problem can be seen as a generalization of the Oberwolfach problem, in \cite{Costa20} the problem $FP(F)$ was denoted by $OP(F)$.}
\end{Problem}

As far as we know, there are only four papers dealing with the Factorization Problem for infinite complete graphs, and two of them, concern classic designs. 
In \cite{Ko77} it is shown
that there exists a resolvable 2-design whenever $v=|\mathbb{N}|$ and $k$ is finite; these designs have, in addition, a cyclic automorphism group $G$ acting sharply transitive on the vertex set; briefly they are $G$-regular.
In \cite{DaHoWe14} it is shown that every infinite $2$-design with $k<v$ is necessarily resolvable, and when $k=v$, both resolvable and non-resolvable designs exist.
We point out that both these papers deal with $t$-designs with $t\geq 2$. 

Furthermore, 
in \cite{BoMa10} the authors construct a $G$-regular 1-factorization of a countable complete graph
for every finitely generated abelian infinite group $G$.
Finally, \cite{Costa20} proves the following.
\begin{Theorem}\label{Costa} Let $F$ be a graph whose order is the cardinal number $\aleph$.
$FP(F)$ has a $G$-regular solution whenever the following two conditions hold:
\begin{enumerate} 
  \item\label{locallyfinite} $F$ is locally finite,
  \item $G$ is an involution free group of order $\aleph$.
\end{enumerate}
\end{Theorem}
Note that this result generalizes the one obtained in \cite{Ko77} to any complete graph of infinite order $\aleph$, blocks of any size less than $\aleph$, and groups $G$ not necessarily cyclic.
Furthermore, Theorem \ref{Costa} can also be seen as a generalization of the result in \cite{BoMa10} to complete graphs of any infinite order.\\

When dealing with infinite graphs, a central role is played by the Rado graph $R$
(see \cite{Rado}), 
named after Richard Rado who gave one of its first explicit constructions. Indeed,
$R$ is the unique countably infinite random graph, and it can be constructed as follows:
$V(R)=\mathbb{N}$ and a pair $\{i,j\}$ with $i<j$ is an edge of $R$ if and only if 
the $i$-th bit of the binary representation of $j$ is one. $R$ shows many interesting  properties, such as the universal property: every finite or countable graph can be embedded as an induced subgraph of $R$. 

When replacing the concept of induced subgraph with the dual one of factor, a weaker result holds. Indeed, in \cite{Ca97} it is pointed out that a countable graph $F$
can be embedded as a factor of $R$ if and only if the domination number of $F$ is infinite. In the same paper, it is further shown that $FP(\mathcal{F}, R)$ has a solution whenever $\mathcal{F}$ is infinite and each of its graphs is locally finite.
Note that a locally finite countable graph has infinite domination number, 
but the converse is not true: 
for example, the Rado graph is not locally finite and it has no finite dominating set
(indeed, for every $D=\{i_1, \ldots, i_t\}\subset\mathbb{N}$, there exists an integer $j\in\mathbb{N}$ whose binary representation has 0 in positions $i_1, \ldots, i_t$, which means  that $j$ is adjacent with no vertex of $D$).

In this paper, we extend this result to any countable family $\mathcal{F}$ of admissible graphs. More precisely, we prove the following. 

\begin{Theorem}\label{FactorizationRado}
Let $\mathcal{F}$ be a countable family of countable graphs. Then, 
$FP(\mathcal{F},R)$ has a solution if and only if the domination number of each graph of $\mathcal{F}$ is infinite.
\end{Theorem}

Furthermore, we prove the solvability of $FP(\mathcal{F})$ whenever the size of 
$\mathcal{F}$ coincides with the order and the domination number of its graphs.

\begin{Theorem}\label{main} 
Let $\mathcal{F}$ be a family of graphs, each of which has order $\aleph$. 
$FP(\mathcal{F})$ has a solution whenever the following two conditions hold:
\begin{enumerate} 
  \item $|\mathcal{F}|=\aleph$, \label{maincondition1}
  \item
  the domination number of each graph in $\mathcal{F}$ is $\aleph$.\label{maincondition2}
\end{enumerate}
\end{Theorem}
When $\mathcal{F}$ contains only copies of a given graph $F$ satisfying condition \ref{locallyfinite} of Theorem \ref{Costa} (i.e., $F$ is locally finite), then $\mathcal{F}$ satisfies both conditions 1 and 2 of Theorem \ref{main}.
Therefore, Theorem \ref{main} can be seen as a generalization of Theorem \ref{Costa}, even though it does not provide any information on the automorphisms of a solution to FP. 
 
Note that if we just require that the domination number of each graph of $\mathcal{F}$ is $\aleph$, there may exist factorizations with fewer factors than $\aleph$; this means that the two conditions in Theorem \ref{main} are independent.
Indeed, the Rado graph $R$ has no finite dominating set
and Corollary \ref{CopiesRado} shows that for every $n\geq 2$ there exists a factorization of $K_\mathbb{N}$ into $n$ copies of $R$. 
We point out that Theorem \ref{main} constructs instead  factorizations of $K_\mathbb{N}$ into infinite copies of $R$.



The paper is organized as follows. 
In Sections 2 and 3, we prove the main results of this paper, Theorems \ref{FactorizationRado} and \ref{main}. 
In Section 4, we deal with $F$-factorizations of $K_{\mathbb{N}}$ when $F$ belongs to a special class of graphs with finite domination number (and hence not satisfying condition \ref{maincondition2} of Theorem \ref{main}): the countable $k$-stars 
(briefly, $S_k$), that is, the vertex disjoint union of $k$ countable stars.
We prove that $FP(S_k)$ has a solution whenever $k>3$, and there is no solution for $k\in \{1,2\}$.
This shows that there are families $\mathcal{F}$ of graphs for which 
$FP(\mathcal{F})$ is not solvable. We leave open the problem when $k=3$. 
In the last section, inspired by \cite{DaHoWe14}, we provide a sufficient condition for a decomposition $\mathcal{F}$ of $K_{\aleph}$ to be resolvable (i.e., the graphs of 
$\mathcal{F}$ can be partitioned into factors of $K_{\aleph}$).


\section{Factorizing the Rado graph}\label{factorizingR}
In this section, we prove Theorem \ref{FactorizationRado}. Also, since the Rado graph $R$ is self-complementary, that is, 
$K_\mathbb{N}\setminus R$ is isomorphic to $R$, we obtain as a corollary 
the countable version of Theorem \ref{main}.

We start by recalling an important characterization of the Rado graph
(see, for example, \cite{Ca97}).

\begin{Theorem}\label{Fact2}
A countable graph is isomorphic to the Rado graph if and only if it satisfies the following property:
\begin{itemize}
\item[$\star$] for every disjoint finite sets of vertices $U$ and $W$, there exists a vertex $z$ adjacent to all the vertices of $U$ and non-adjacent to all the vertices of $V$.
\end{itemize} 
\end{Theorem}

Now we slightly generalize the construction of the Rado graph given in the introduction.
\begin{Definition}\label{RqI:defi}
Given a set $I\subset\{0,\dots,q-1\}$, 
with $1 \leq |I| < q$, we denote by $R^q_{I}$ the following graph:
$V(R^q_{I})=\mathbb{N}$, and $\{x,y\}$, with $x<y$, is an edge of $R^q_{I}$ whenever the $x$-th digit of $y$ in the base $q$ expansion of $y$ belongs to $I$.  
\end{Definition}

Cleraly, when $q=2$ and $I=\{1\}$
we obtain our initial description of the Rado graph (i.e. $R=R^2_{\{1\}}$).

\begin{Proposition}\label{RqI}
Every graph $R^q_{I}$ is isomorphic to the Rado graph. 
\end{Proposition}
\begin{proof} By Theorem \ref{Fact2}, it is enough to show 
that property $\star$ holds for $R^q_{I}$. 
We assume, without loss of generality, that 
$0\in I$ while $1\not\in I$, and let $U$ and $V$ be two disjoint subsets of $\mathbb{N}$. Then there are infinitely many positive integers whose base $q$ expansion has $0$ in each position $u\in U$
and $1$ in each position $v\in V$. Denoting by $z$ one of these integers larger than $\max(U \cup V)$, we have that 
$z$ is adjacent to all the vertices of $U$ but to none in $V$.
\end{proof}

Note that $K_{\mathbb{N}}= \bigcup_{i=0}^{q-1} R^q_{\{i\}}$
and $R^q_{\{0,\ldots,q-2\}}=\bigcup_{i=0}^{q-2} R^q_{i}$.
Considering that the $R^q_{\{i\}}$s are pairwise edge-disjoint and isomorphic to the Rado graph, we obtain the following.

\begin{Corollary}\label{CopiesRado}
For every $n\in\mathbb{N}$, the graphs $R$ and $K_{\mathbb{N}}$ can be factorized
into $n$ and $n+1$ copies of $R$, respectively.
\end{Corollary}

The following result is crucial to prove Theorem \ref{FactorizationRado}. It strengthens a result given in \cite{Ca97} and allows us to suitably embed in the Rado graph $R$ 
any countable graph with infinite domination number.

\begin{Proposition}\label{embedding}
Let $F$ be a countable graph with no finite dominating set. 
For every edge $e\in E(R)$, there exists an embedding $\sigma_e$ of $F$ in $R$ such that:
\begin{enumerate}
\item\label{embedding1} $\sigma_e(F)$ is a spanning subgraph of $R$ containing the edge $e$;
\item\label{embedding2} $R\setminus \sigma_e(F)$ is isomorphic to $R$.
\end{enumerate}
\end{Proposition}
\begin{proof} By Proposition \ref{RqI}, the graphs $R_{\{0,1\}}^3$,
$R_{\{0\}}^3$ and $R_{\{1\}}^3$ are isomoprhic to $R$. 
Therefore, we can take $R= R_{\{0,1\}}^3$. 

Let $e$ be an edge of $R  = R_{\{0\}}^3 \cup R_{\{1\}}^3$. 
We can assume without loss of generality that 
$e$ lies in $R_{\{0\}}^3$.
In  \cite[Proposition 8]{Ca97}, 
it is shown that there exists an embedding 
$\sigma_e$ of $F$ into the Rado graph $R_{\{0\}}^3 \subset R$
satisfying condition $1$. It is then left to 
prove that condition \ref{embedding2} holds. By Theorem \ref{Fact2}, this is equivalent to saying that 
$R\setminus \sigma_e(F)$ satisfies $\star$. 

Let $U$ and $V$ be two finite disjoint subsets of $\mathbb{N}$. 
Clearly, there are infinitely many positive integers whose base $3$ expansion has $1$ in each position $u\in U$ and $2$ in each position $v\in V$.
Let $z$ be one of these integers larger than $\max(U \cup V)$.
Since $R\setminus \sigma_e(F)$ contains $R_{\{1\}}^3$ and it is edge-disjoint with $R_{\{2\}}^3$, 
it follows that  $z$ is adjacent in $R\setminus \sigma_e(F)$ to all the vertices of $U$ and is non-adjacent to all the vertices of $V$.
This means that $\star$ holds for $R\setminus \sigma_e(F)$.
\end{proof}

We are now ready to prove Theorem \ref{FactorizationRado}, whose statement is recalled here, for clarity.\\

\noindent
\textbf{Theorem \ref{FactorizationRado}.}
\emph{Let $\mathcal{F}$ be a countable family of countable graphs. Then, $FP(\mathcal{F},R)$ has a solution if and only if the domination number of each graph of $\mathcal{F}$ is infinite.}
\begin{proof}
Since the Rado graph has no finite dominating set, the same holds for its spanning subgraphs. Hence, each graph of $\mathcal{F}$ must have infinite domination number. Under this assumption, we are going to show that $FP(\mathcal{F},R)$ has a solution.
 
Let $E(R)=\{ e_1,\dots,e_n,\ldots\}$ and 
$\mathcal{F}=\{F_1,\dots,F_n,\ldots\}$. By recursively applying Proposition \ref{embedding}, we obtain a sequence
of isomorphisms $\sigma_i: F_i\rightarrow \Gamma_i$ satisfying 
for each $i\in\N$ the following properties:
\begin{itemize}
\item $\Gamma_i$ is a spanning subgraph of $R$;
\item $R\setminus (\Gamma_1\cup\Gamma_2\cup \cdots\cup \Gamma_{i-1})$
 is isomorphic to $R$ and contains $\Gamma_i$;
\item $e_i$ lies in $\Gamma_1\cup\Gamma_2\cup \cdots\cup \Gamma_i$.
\end{itemize}
It follows that the $\Gamma_i$s are pairwise edge-disjoint factors of $R$ which partition $E(R)$. Therefore, 
$\{\Gamma_i: i\in \mathbb{N}\}$ is a solution to $FP(\mathcal{F},R)$.
\end{proof}

The proof of Theorem \ref{FactorizationRado} allows us to construct solutions to $FP(\mathcal{F},R)$ even when
the cardinality of $\mathcal{F}$ is finite, 
provided that $\mathcal{F}$ contains a copy of the Rado graph. In other words, we have the following.

\begin{Corollary}\label{FactorizationRado:coro}
Let $\mathcal{F}$ be a finite family of countable graphs such that
\begin{enumerate} 
  \item $\mathcal{F}$ contains at least one graph isomorphic to the Rado graph;
  \item  the domination number of each graph in $\mathcal{F}$ is infinite.
\end{enumerate}
Then, $FP(\mathcal{F},R)$ has a solution.
\end{Corollary}

Recalling that $R$ is self complementary, the countable version of Theorem \ref{main} can be easily obtained as a corollary to Theorem \ref{FactorizationRado}.
\begin{Corollary}\label{FactorizationRado:coro:2}
Let $\mathcal{F}$ be a countable family of countable graphs. 
$FP(\mathcal{F})$ has a solution whenever the domination number of each graph in $\mathcal{F}$ is infinite.
\end{Corollary}
\begin{proof}
Recall that $R_{\{0\}}^2$ and 
$R_{\{1\}}^2$ are copies of $R$ which together factorize $K_{\mathbb{N}}$.   Therefore, it is enough to
partition $\mathcal{F}$ into two countable families $\mathcal{F}_1$ and $\mathcal{F}_2$, and then apply 
Theorem \ref{FactorizationRado} to get a solution 
$\mathcal{G}_i$ to $FP(\mathcal{F}_i, R_{\{i\}}^2)$, for $i=0,1$.
Clearly, $\mathcal{G}_1\cup \mathcal{G}_2$ provides a solution to  $FP(\mathcal{F})$.
\end{proof}

The natural generalization of property $\star$ to a generic cardinality $\aleph$ is the following one:

\begin{itemize}
\item[$\star_{\aleph}$] for every disjoint sets of vertices $U$ and $W$ whose cardinality is smaller than $\aleph$, there exists a vertex $z$ adjacent to all the vertices of $U$ and non-adjacent to all the vertices of $V$.
\end{itemize}
 Then, using the transfinite induction (see Theorem \ref{transfinite} below), one could also prove the following generalization of Proposition \ref{Fact2}:
 
\begin{Proposition}
Any two graphs of order $\aleph$ that satisfy property $\star_{\aleph}$ are pairwise isomorphic. 
\end{Proposition}

Therefore, we can refer to any graph of order $\aleph$ 
and satisfying property $\star_{\aleph}$ 
as the  $\aleph$-Rado graph $R_{\aleph}$. 
Its existence is guaranteed under the 
Generalized Continuum Hypothesis (GCH) which states that if 
$\aleph'\prec \aleph$ then $2^{\aleph'}\preceq \aleph$.
Under GCH, one can see that the set of all $q$-ary sequences of length $\prec\aleph$ 
has size $\aleph$:  indeed, for every  $\aleph'\prec \aleph$, the set of all 
$q$-ary sequences of length $\aleph'$ has cardinality $2^{\aleph'}$, and by GCH
we have that $2^{\aleph'}\preceq \aleph$. 
This means that the construction of the countable Rado graph 
(Definition \ref{RqI:defi}) based on representing every natural number with a finite $q$-ary sequence (its base $q$ expansion) can be generalized to any order.

By assuming that GCH holds, we can prove the following generalization of Theorem \ref{FactorizationRado}.

\begin{Theorem}\label{FactorizationRadoAleph}

Let $\mathcal{F}$ be a family of graphs of order $\aleph$  
and assume that $|\mathcal{F}| = \aleph$. 
Then $FP(\mathcal{F},R_{\aleph})$ has a solution if and only if 
the domination number of each graph in $\mathcal{F}$ is $\aleph$.
\end{Theorem}


\section{Factorizing infinite complete graphs}
\label{factorizingK}
In this section we prove Theorem \ref{main}.
We point out that if we assume the Generalized Continuum Hypothesis, 
considering that $R_\aleph$ is self complementary by property $\star_\aleph$,
we can proceed as in Corollary \ref{FactorizationRado:coro:2}
and obtain Theorem \ref{main} as a consequence of
Theorem \ref{FactorizationRadoAleph}.

Here we present a proof of Theorem \ref{main}
that does not require GCH, which we recall is independent of ZFC. 
Therefore, a proof that does not require GCH is to be preferred.

We say that a graph or a set of vertices is $\aleph$-small (resp. $\aleph$-bounded) if their order or cardinality is smaller than $\aleph$ (resp. smaller than or equal to $\aleph$).
Given two graphs $F$ and $\Lambda$ of order $\aleph$, we denote by $\Sigma_{\aleph}(F, \Lambda)$ the set of all graph embeddings 
between an induced $\aleph$-small subgraph of $F$ and a subgraph of $\Lambda$.
A partial order on $\Sigma_{\aleph}(F, \Lambda)$ can be easily defined as follows:
if 
$\sigma: G \rightarrow\Gamma$ and 
$\sigma':G'\rightarrow\Gamma'$ are embeddings of $\Sigma_{\aleph}(F, \Lambda)$, we say that
$\sigma\leq \sigma'$ whenever $\sigma'$ is an extension of $\sigma$, namely, 
$G$ and $\Gamma$ are subgraphs of  $G'$ and $\Gamma'$, respectively,
and $\sigma'|_{G} = \sigma$ (where $\sigma'|_{G}$ is the restriction of $\sigma'$ to $G$).

\begin{Lemma}\label{Estensione1}
Let $F$ be a graph of order $\aleph$ and with no $\aleph$-small dominating set. 
Also, let $\Theta$ be an $\aleph$-small subgraph of $K_{\aleph}$,  
and let $\sigma\in \Sigma_{\aleph}(F, K_{\aleph}\setminus\Theta)$.
\begin{enumerate}
\item
If $v \in V(F)$, then there is an embedding $\sigma': G' \rightarrow \Gamma'$ in 
$\Sigma_{\aleph}(F, K_{\aleph}\setminus\Theta)$ such that 
\[|V(G')|\leq |V(G)|+1, \;\;\sigma\leq \sigma'\;\;\text{and}\;\; v\in V(G');\]
\item
If $x \in V(K_{\aleph})$, then there is an embedding $\sigma'': G'' \rightarrow \Gamma''$ in 
$\Sigma_{\aleph}(F, K_{\aleph}\setminus\Theta)$ such that
\[|V(G'')|\leq|V(G)|+1, \;\; \sigma\leq \sigma''\;\;\text{and}\;\; x\in V(\Gamma'').\]
\end{enumerate}
\end{Lemma}
\begin{proof}
Let $\sigma: G \rightarrow \Gamma$ be an embedding in
$\Sigma_{\aleph}(F, K_\mathbb{N}\setminus\Theta)$, and let $v\in V(F)$ and $x\in V(K_{\aleph})$.
Clearly, when $v\in V(G)$ or $x\in V(\Gamma)$, 
we can take $\sigma'=\sigma$ or $\sigma''=\sigma$, respectively. Therefore, we can assume $v\not \in V(G)$ and $x\not\in V(\Gamma)$.
\begin{enumerate} 
\item 
Let $G'$ be the subgraph of $F$ induced by $v$ and $V(G)$. 
Since $V(\Theta)$ is $\aleph$-small, we can  choose $a\in V(K_{\aleph})\setminus V(\Theta)$ and
let $\sigma': V(G)\ \cup \{v\} \rightarrow V(\Gamma)\ \cup\ \{a\}$ be the extension of
$\sigma$ such that $\sigma'(v)=a$.  Setting $\Gamma' = \sigma'(G')$, we have that
$\sigma'$ is the required embedding of $\Sigma_{\aleph}(F, K_{\aleph}\setminus\Theta)$.

\item
Since $F$ has no $\aleph$-small dominating set, $V(G)$ (which is an $\aleph$-small set) 
cannot be a dominating set for $F$. 
Hence, there is a vertex $a\in V(F)$ that is not adjacent to any of the vertices
of $G$. We denote by $G''$ (resp., $\Gamma''$) the graph obtained by adding 
$a$ to $G$ (resp., $x$ to $\Gamma)$ as an isolated vertex. Clearly, $G''$ is an induced subgraph of $F$; also, $\Gamma''$ and $\Theta$ have no edge in common, 
since $E(\Gamma'')= E(\Gamma)$. 
Therefore, the extension $\sigma'': G'' \rightarrow \Gamma''$ of $\sigma$
such that $\sigma''(a)=x$ is the required embedding of 
$\Sigma_{\aleph}(F, K_\aleph\setminus\Theta)$.
\end{enumerate}
\end{proof}

From now on, we will work within the Zermelo-Frankel axiomatic system with the Axiom of Choice in the form of the Well-Ordering Theorem. We recall the definition of a well-order.

\begin{Definition}
A well-order $\prec$ on a set $X$ is a total order on $X$ with the property that every non-empty subset of $X$ has a least element.
\end{Definition}

The following theorem is equivalent to the Axiom of Choice.

\begin{Theorem}[Well-Ordering]
Every set $X$ admits a well-order $\prec$.
\end{Theorem}

Given an element $x\in X$, we define the section $X_{\prec x}$ associated to it:
$$X_{\prec x}=\{y\in X: y\prec x\}.$$

\begin{Corollary}\label{goodgoodorder}
Every set $X$ admits a well-order $\prec$ such that the cardinality of any section is smaller than $|X|$.
\end{Corollary}
\begin{proof}
Let us consider a well-order $\prec$ on $X$. Let $x$ be the smallest element such that $X_{\prec x}$ has the same cardinality as $X$. The set $Y=X_{\prec x}$ is such that all its sections with respect to the order $\prec$ have smaller cardinality. Since $Y$ instead has the same cardinality as $X$, the order $\prec$ on $Y$ induces an order $\prec'$ on $X$ with the required property.
\end{proof}

We recall now that well-orderings allow proofs by induction.

\begin{Theorem}[Transfinite induction]\label{transfinite}
Let $X$ be a set with a well-order $\prec$ and let $P_x$ denote a property for each 
$x\in X$. Set $0=\min  X$ and assume that:
\begin{itemize}
\item $P_0$ is true, and
\item for every $x\in X$, if $P_y$ holds for every $y\in X_{\prec x}$, then $P_x$ holds.
\end{itemize}
Then $P_x$ is true for every $x\in X$.
\end{Theorem}

We are now ready to prove Theorem \ref{main}. 
The idea behind the proof can be better understood by restricting our attention to the countable case, $\aleph = \mathbb{N}$. 
To solve $FP(\{F_\alpha: \alpha\in \mathbb{N}\})$, 
we first order the edges of $K_\mathbb{N}: \{e_0, e_1, \ldots, e_\gamma, \ldots\}$. 
Then, we define
embeddings $\sigma_\alpha^\beta: G_\alpha^\beta\rightarrow \Gamma_\alpha^\beta$
where $G_\alpha^\beta$ is an induced subgraph of $F_\alpha$, 
and $\Gamma_\alpha^\beta$ is a subgraph of $K_{\mathbb{N}}$. 
These embeddings are obtained by recursively applying Lemma \ref{Estensione1}
which adds, at each step, a vertex to  $G_\alpha^\beta$ and  a vertex to  $\Gamma_\alpha^\beta$
and makes sure that the vertex $\beta$ belongs to both these graphs
(this procedure can be seen as a variation of Cantor's ``back-and-forth'' method).
We also make sure that, for every $\gamma$, 
the graphs $\Gamma_0^\gamma, \Gamma_1^\gamma, \ldots, \Gamma_\gamma^\gamma$ are pairwise edge-disjoint
and contain between them the edge $e_\gamma$. 
The solution to $FP(\{F_\alpha: \alpha\in \mathbb{N}\})$ will be represented by 
$\mathcal{G}=\{\Gamma_\alpha: \alpha\in \mathbb{N}\}$ where $\Gamma_\alpha=\bigcup_{\beta} \Gamma_\alpha^\beta$.
 \\

\noindent
\textbf{Theorem \ref{main}.}
\emph{
Let $\mathcal{F}$ be a family of graphs, each of which has order $\aleph$. 
$FP(\mathcal{F})$ has a solution whenever the following two conditions hold:
\begin{enumerate} 
  \item $|\mathcal{F}|=\aleph$, 
  \item
  the domination number of each graph in $\mathcal{F}$ is $\aleph$.
\end{enumerate}
}
\begin{proof}
Let $\mathcal{F} = \{F_\alpha:\alpha\in\mathcal{A}\}$.
We consider a well-order $\prec$ on $\mathcal{A}$ satisfying Corollary \ref{goodgoodorder}.
Since by assumption $|V(F_{\alpha})|=|\mathcal{A}|=\aleph$, for every $\alpha\in \mathcal{A}$, we can take $V(F_{\alpha})=V(K_{\aleph})= \mathcal{A}$ and index the edges of $K_{\aleph}$ over $\mathcal{A}$: $E(K_{\aleph})=\{e_{\alpha}: \alpha\in  \mathcal{A}\}$.

To prove the assertion, we construct a chain of families 
$(\mathcal{E}_{\gamma})_{\gamma\in \mathcal{A}}$, where 
\[
\mathcal{E}_{\gamma}:=\{\sigma_\alpha^\beta: G_\alpha^{\beta} \rightarrow \Gamma_\alpha^\beta\mid  \sigma_\alpha^\beta \in \Sigma_{\aleph}(F_\alpha, K_{\aleph}),(\alpha,\beta) \in \mathcal{A}_{\preceq \gamma}\times\mathcal{A}_{\preceq \gamma}\},
\]
which satisfy the ascending property, that is, 
$\mathcal{E}_{\gamma'}\subseteq \mathcal{E}_{\gamma}$ if $\gamma'\preceq \gamma$, 
and the following three conditions: 
  \begin{enumerate}
    \item[($1_\gamma$)] for every $(\alpha,\beta)\in\mathcal{A}_{\preceq \gamma}\times\mathcal{A}_{\preceq \gamma}$ and $\beta'\prec \beta$ we have that $\sigma_{\alpha}^{\beta'}\leq \sigma_{\alpha}^{\beta}$ and $\beta\in V(G^{\beta}_{\alpha})\ \cap \ V(\Gamma^{\beta}_{\alpha})$;
    \item[($2_\gamma$)] for every $\beta\in\mathcal{A}_{\preceq \gamma}$, the graphs 
    $\Gamma_{\alpha}^\beta: \alpha\preceq \beta $ are pairwise edge-disjoint, and the edge $e_\beta$ belongs to their union;
    \item[($3_\gamma$)] for every $\alpha,\beta\in\mathcal{A}_{\preceq \gamma}$, the graph 
    $\Gamma_{\alpha}^\beta$ is either finite or $|\mathcal{A}_{\preceq \gamma}|$-bounded.
  \end{enumerate}
The desired factorization of $K_\aleph$ is then
$\mathcal{G}=\{\Gamma_\alpha: \alpha\in \mathcal{A}\}$, 
where $\Gamma_\alpha =\bigcup_{\beta\in\mathcal{A}} \Gamma^\beta_\alpha$ for every $\alpha\in\mathcal{A}$.
Indeed, properties ($1_\gamma$) guarantee that each $\Gamma_\alpha$ is a factor of $K_\aleph$ isomorphic to $F_\alpha$. Also, properties ($2_\gamma$) ensure that the $\Gamma_\alpha$s are pairwise edge-disjoint and between them 
contain all the edges of $K_{\aleph}$.

We proceed by transfinite induction on $\gamma$.

BASE CASE. Let $0=\min  \mathcal{A}$, 
choose an edge $e\in E(F_{0})$ and let $\sigma\in \Sigma_\aleph(F_{0}, K_\aleph)$ 
be the embedding that maps $e$ to $e_{0}$. By Lemma \ref{Estensione1}, there exists an embedding 
$\sigma_0^0: G_0^0 \rightarrow \Gamma_0^0$ in 
$\Sigma_\aleph(F_{0}, K_\aleph)$ such that $\Gamma_0^0$ is a finite graph and
\[ \sigma\leq \sigma_0^0\;\;\text{and}\;\; 0\in 
V(G_0^0)\cap V(\Gamma_0^0).
\] 
Clearly, $\mathcal{E}_0:=\{\sigma_{0}^0\}$ satisfies properties ($1_0$), ($2_0$) and ($3_0$).\\

TRANSFINITE INDUCTIVE STEP.
We assume that, for any $\gamma'\prec \gamma$, there is a family $\mathcal{E}_{\gamma'}$ satisfying properties ($1_{\gamma'}$), ($2_{\gamma'}$) and ($3_{\gamma'}$), and
prove that it can be extended to a family $\mathcal{E}_{\gamma}$ that satisfies properties ($1_{\gamma}$), ($2_{\gamma}$) and ($3_{\gamma}$).
Clearly it is enough to provide the maps $\sigma_\alpha^\beta$ where either $\alpha=\gamma$ or $\beta=\gamma$. 

We start by constructing the maps $\sigma_{\alpha}^{\gamma}$ for every $\alpha\prec \gamma$. 
We proceed by transfinite induction on $\alpha$.
\begin{itemize}
\item Base case. Set $\Theta_0:=\bigcup_{\alpha, \beta\prec \gamma} \Gamma_\alpha^\beta$ and note that, by property ($3_{\gamma'}$), $\Theta_0$ is $\aleph$-small. We also set $\sigma_{0}^{\prec\gamma}: \bigcup_{\beta\prec \gamma} G_0^{\beta} \rightarrow \bigcup_{\beta\prec \gamma}  \Gamma_0^{\beta}$ to be the map of 
$\Sigma_{\aleph}(F_0, K_\aleph\setminus\Theta_0)$ whose restriction to $G_0^{\beta}$ is $\sigma_0^{\beta}$. We note that property ($3_{\gamma'}$) guarantees that the order of $\bigcup_{\beta\prec \gamma} G_0^{\beta}$ is either finite or $|\mathcal{A}_{\preceq \gamma}|$-bounded, 
hence $\aleph$-small. 

Therefore, we can apply Lemma \ref{Estensione1} 
(with  $\sigma = \sigma_{0}^{\prec\gamma}$) to obtain 
the map $\sigma_0^{\gamma}: G_0^{\gamma} \rightarrow \Gamma_0^{\gamma}$ in 
$\Sigma_{\aleph}(F_0, K_\aleph\setminus\Theta_0)$ such that $|V(\Gamma_0^{\gamma})|\leq|V(\bigcup_{\beta\prec \gamma} \Gamma_0^{\beta})|+2$ and, for every $\gamma'\prec \gamma$,
\[ \sigma_0^{\gamma'}\leq \sigma_0^{\gamma}\;\;\text{and}\;\; \gamma\in V(G_0^{\gamma})\cap V(\Gamma_0^{\gamma}).\] 
\item Inductive step. Assume we have defined the maps 
$\sigma_{\alpha'}^{\gamma}$ for every $\alpha'\prec \alpha$, and set 
\[
\Theta_\alpha:=\bigcup_{\alpha'\prec \alpha} \Gamma_{\alpha'}^{\gamma}\;\cup\; 
          \bigcup_{\alpha \prec\alpha'\prec \gamma, \beta \prec \gamma} \Gamma_{\alpha'}^{\beta}.
\] 
As before,  by Lemma \ref{Estensione1} there exists $\sigma_\alpha^{\gamma}: G_\alpha^{\gamma} \rightarrow \Gamma_\alpha^{\gamma}$ in 
$\Sigma_{\aleph}(F_\alpha, K_\aleph\setminus\Theta_\alpha)$ such that $|V(\Gamma_\alpha^{\gamma})|\leq|V(\bigcup_{\beta\prec \gamma} \Gamma_\alpha^{\beta})|+2$ and, for every $\gamma'\prec \gamma$,
\[ 
\sigma_\alpha^{\gamma'}\leq \sigma_\alpha^{\gamma}\;\;\text{and}\;\;\gamma\in V(G_\alpha^{\gamma})\cap V(\Gamma_\alpha^{\gamma}).
\] 
\end{itemize}
Finally, we define the maps $\sigma_{\gamma}^\beta$ when $\beta \preceq \gamma$. We set $\Theta:=\bigcup_{\alpha \prec \gamma} \Gamma_\alpha^{\gamma}$ and proceed by transfinite induction on $\beta$.
\begin{itemize}
\item Base case. 
If $e_{\gamma}\in \Theta$, let $\sigma$ be the empty map of $\Sigma_{\aleph}(F_{\gamma}, K_\aleph\setminus\Theta)$. Otherwise, chose an edge $e\in E(F_{\gamma})$, 
and let $\sigma\in \Sigma_{\aleph}(F_{\gamma}, K_\aleph\setminus\Theta)$ 
be the embedding that maps $e$ to $e_{\gamma}$. 
By Lemma \ref{Estensione1}, there exists $\sigma_{\gamma}^0: G_{\gamma}^0 \rightarrow \Gamma_{\gamma}^0$ in 
$\Sigma_{\aleph}(F_{\gamma}, K_\aleph\setminus\Theta)$ such that $\Gamma_{\gamma}^0$ is a finite graph and
\[ \sigma\leq \sigma_{\gamma}^0\;\;\text{and}\;\; 0\in V(G_{\gamma}^0)\cap V(\Gamma_{\gamma}^0).\]
\item Inductive step.
Assume we have defined the maps $\sigma_{\gamma}^{\beta'}$ for any $\beta'\prec \beta$. Again by Lemma \ref{Estensione1}, there exists $\sigma_{\gamma}^{\beta}: G_{\gamma}^{\beta} \rightarrow \Gamma_{\gamma}^{\beta}$ in 
$\Sigma_{\aleph}(F_{\gamma}, K_\aleph\setminus\Theta)$ such that $|V(\Gamma_\gamma^{\beta})|\leq|V(\bigcup_{\beta'\prec \beta} \Gamma_\gamma^{\beta'})|+2$ and, for any $\beta'\prec \beta$,
\[ \sigma_{\gamma}^{\beta'}\leq \sigma_{\gamma}^{\beta}\;\;\text{and}\;\; 
\beta\in  
V(G_{\gamma}^{\beta})\cap V(\Gamma_{\gamma}^{\beta}).
\] 
\end{itemize}
It follows from the construction that the family \[\mathcal{E}_{\gamma}:=\{\sigma_\alpha^\beta: G^\alpha_{\beta} \rightarrow \Gamma_\alpha^\beta\mid  \sigma_\alpha^\beta \in \Sigma_\aleph(F_\alpha, K_\aleph), \alpha,\beta\preceq \gamma\}\]
satisfies properties $(1_\gamma)$, $(2_\gamma)$ and $(3_\gamma)$.
\end{proof}

\section{The Factorization Problem for $k$-stars}
Theorem \ref{main} does not provide solutions to $FP(F)$ whenever the graph $F$ has a dominating set 
of cardinality less than its order. In particular, if $F$ is countable with a finite dominating set, then the existence of a solution to $FP(F)$ is an open problem. In this section, we consider a special class of such graphs, the \emph{$k$-stars} $S_k$. More precisely, 
\begin{itemize}
\item the star $S_1$ is the graph with vertex-set $\mathbb{N}$ whose edges are of the form 
$\{0, i\}$ for every $i\in \mathbb{N}\setminus \{0\}$; 
\item the $k$-star $S_{k}$ is the vertex-disjoint union of $k$ stars. 
\end{itemize}
Note that $S_k$ contains exactly $k$ vertices of infinite degree, which we call \emph{centers} 
and form a finite dominating set of $S_k$. 

In the following, we show that $FP(S_k)$ has no solution whenever $k\in \{1,2\}$, while it admits a solution for every $k>3$. Unfortunately, we leave open the problem for $3$-stars.  

\subsection{The case $k\in \{1,2\}$}
\begin{Proposition}\label{k1}
$FP(S_1)$ has no solution.
\end{Proposition}
\begin{proof}
Assume for a contradiction that there is a factorization $\mathcal{G}$ of $K_{\mathbb{N}}$ into $1$-stars. Choose any star $\Gamma\in \mathcal{G}$ and let $g$ denote its center. Considering that
all the edges of $K_{\mathbb{N}}$ incident with $g$ belong $\Gamma$, it follows that $g$ cannot be a vertex in any other star of $\mathcal{G}$, which therefore are not factors and this is a contradiction.
\end{proof}

With essentially the same proof, one obtains the following.
\begin{Remark}
Let $F$ be the vertex-disjoint union of $S_1$ with a finite set of isolated vertices. 
Then $FP(F)$ has no solution.
\end{Remark}

To prove the non-existence of a solution to $FP(S_2)$ it will be useful the following lemma.

\begin{Lemma}\label{everyvertexiscenter}
If $\mathcal{G}$ is a factorization of $K_{\mathbb{N}}$ into $k$-stars, then
there is at most one vertex of $K_{\mathbb{N}}$ that is never a center in any $k$-star of $\mathcal{G}$. It follows that $|\mathcal{G}|=|\mathbb{N}|$.
\end{Lemma}
\begin{proof} It is enough to notice that every pair $\{a,b\}$ of vertices of 
$K_{\mathbb{N}}$ is the edge of some $2$-star $\Gamma$ of $\mathcal{G}$; 
hence, either $a$ or $b$ is a center of $\Gamma$.
\end{proof}

\begin{Proposition}\label{k2}
$FP(S_2)$ has no solution.
\end{Proposition}
\begin{proof}
Assume for a contradiction that there is a factorization $\mathcal{G}$ of $K_{\mathbb{N}}$ into $2$-stars. For every $\Gamma\in\mathcal{G}$, letting $c$ be a center of $\Gamma$, we denote by 
$\Gamma(c)$ the set of vertices adjacent with $c$ in $\Gamma$ (i.e., the neighborhood of $c$ in $\Gamma$).

Choose any $2$-star $\Gamma\in \mathcal{G}$ and let $a$ and $b$ denote its centers. 
Also, let $\Gamma'$ be the $2$-star of $\mathcal{G}\setminus\{\Gamma\}$ containing the edge $\{a,b\}$. 
Without loss of generality, we can assume that $a$ is a center of $\Gamma'$. 
Finally, by Lemma \ref{everyvertexiscenter}, we can choose $x\in \Gamma'(a)\setminus\{b\}$ 
such that there exists a $2$-star  $\Gamma''\in \mathcal{G}$ having $x$ as a vertex.

Since $\Gamma$ is a factor of 
$K_{\mathbb{N}}$, it follows that $x\in \Gamma(b)$. In other words, 
$\Gamma\ \cup\ \Gamma'$ contains the edges $\{x,a\}$ and $\{x,b\}$. Therefore,
$a,b\not\in\Gamma''(x)$. Since $\Gamma''$ is a factor of $K_{\mathbb{N}}$ 
and $\{a,b\}$ is an edge of $\Gamma$, it follows 
$a,b\in \Gamma''(y)$, where $y$ is the other center of $\Gamma''$.  
In other words, $\{y,a\}$ and $\{y,b\}$ belong to $\Gamma''$, hence $y$
cannot lie in $\Gamma$, contradicting the fact that $\Gamma$ is a factor.
\end{proof}

\subsection{The case $k\geq 4$}
In this section we prove the solvability of $FP(S_k)$ whenever $k\geq 4$. 
For our constructions we need to introduce the following notation.

Let $\mathbb{D}$ be an integral domain and set $V=\mathbb{D}\times \{0,1,\ldots,h\}$, for $h\geq 0$. 
For the sake of brevity, we will denote each pair $(a,i)\in V$ by $a_i$.
Given a graph $\Gamma$ with vertices in $V$, for every $a,b\in \mathbb{D}$ we denote by 
$a\Gamma+b$ the graph obtained by replacing each vertex $x_i$ of $\Gamma$ with $(ax+b)_i$.
Also, we denote by $Orb_{\mathbb{D}}(\Gamma)=\{\Gamma+d:d\in \mathbb{D}\}$ the 
$\mathbb{D}$-orbit of $\Gamma$, that is, the set of all translates of $\Gamma$ by the elements of 
$\mathbb{D}$.

\begin{Proposition}\label{k>3}
For every $k\geq 4$, there exists a $k$-star $\Gamma$ with vertex set $V=\mathbb{Z}\times\{0,1\}$ such that $Orb_{\mathbb{Z}}(\Gamma)$ is a factorization of $K_{V}$ into $k$-stars.
\end{Proposition}
\begin{proof} We first deal with the case $k=4$. Set $\Gamma= \bigcup_{i=1}^4 \Gamma_i$, where
each $\Gamma_i$ is the $1$-star with vertices in  $V=\mathbb{Z}\times\{0,1\}$
and center $x_i$ defined as follows (see Figure \ref{fig1}): 
\begin{itemize}
\item  $x_1= 0_0$ and $\Gamma_1(x_1)=\{i_0: i\geq 1\}$;
\item  $x_2= -1_1$ and $\Gamma_2(x_2)=\{i_1: i\geq 0\}\ \cup\ \{-1_0\}$;
\item  $x_3= -2_0$ and $\Gamma_2(x_3)=\{i_1: i\leq -3\}$;
\item  $x_4= -2_1$ and $\Gamma_4(x_4)=\{i_0: i\leq -3\}$.
\end{itemize}
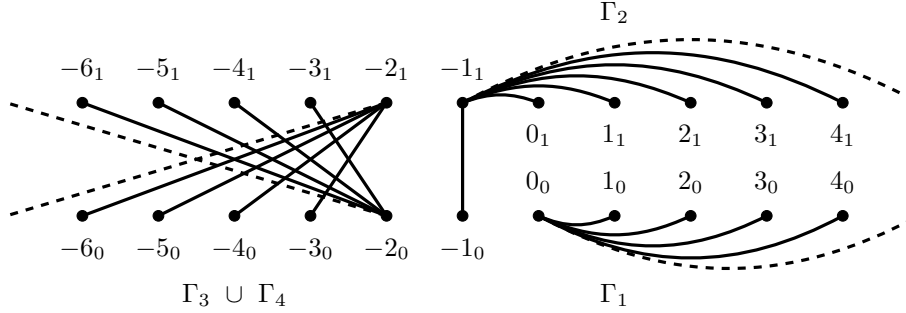
\begin{figure}[h]
\centering
%
\begin{minipage}{15cm}
\hspace{-1.5cm}

\begin{tikzpicture}[x=2cm,y=1cm,scale=0.50]
\foreach \j in {0,1,2,3,4,5,6,7,8,9,10,11,12} 
{
  \coordinate  (a\j) at (0+\j, 3); 
  \coordinate  (b\j) at (0+\j, 0);
}; 
\foreach \j in {1,2,3,4,5} 
{
  \ADD{\j}{-7}{\sol}
  \draw (a\j) circle (4pt)[fill]  node[above=5pt] {$\sol_1$}; 
  \draw (b\j) circle (4pt)[fill]  node[below=5pt] {$\sol_0$}; 
};  
\foreach \j in {7,8,9,10,11} 
{
  \ADD{\j}{-7}{\sol};
  \draw (a\j) circle (4pt)[fill]  node[below=5pt] {$\sol_1$}; 
  \draw (b\j) circle (4pt)[fill]  node[above=5pt] {$\sol_0$}; 
};        
  \draw (a6) circle (4pt)[fill]  node[above=5pt] {$-1_1$}; 
  \draw (b6) circle (4pt)[fill]  node[below=5pt] {$-1_0$}; 
\foreach \j in {1,2,3,4} 
{
  \draw[very thick] (a5) -- (b\j);
  \draw[very thick] (b5) -- (a\j); 
}
  \draw[dashed, very thick] (b5) to (a0);
  \draw[dashed, very thick] (a5) to (b0); 
 
\coordinate  (name1) at (8, -3);   
\draw (name1) node[above=5pt] {$\Gamma_1$};   
\coordinate  (name2) at (8, 4.5);   
\draw (name2) node[above=5pt] {$\Gamma_2$};  
\coordinate  (name3) at (3, -3);   
\draw (name3) node[above=5pt] {$\Gamma_3\ \cup\ \Gamma_4$};

\foreach \j in {7,8,9,10,11} 
{
    \MULTIPLY{1.5}{\j}{\tempA};
    \MULTIPLY{1.5}{7}{\tempB}; 
    \ADD{\tempA}{\tempB}{\solo};       
   \MULTIPLY{-}{\solo}{\tempC}
   \ADD{\tempC}{180}{\soli};
   \draw[very thick, in=\soli, out=\solo]  (a6) to (a\j);
} 
\draw[dashed, very thick, in=151.5, out=28.5]  (a6) to (a12);

\foreach \j in {8,9,10,11} 
{
    \MULTIPLY{1.5}{\j}{\tempA};
    \MULTIPLY{1.5}{8}{\tempB}; 
    \ADD{\tempA}{\tempB}{\solo};       
   \MULTIPLY{-}{\solo}{\tempC}
   \ADD{\tempC}{180}{\soli};
   \draw[very thick, in=-\soli, out=-\solo]  (b7) to (b\j); 
}

   \draw[dashed, very thick, in=-151.5, out=-28.5]  (b7) to (b12);
   \draw[very thick]  (b6) to (a6);   
\end{tikzpicture}

\end{minipage}
\label{fig1}
\end{figure}
We claim that  $\mathcal{G}:=Orb_{\mathbb{Z}}(\Gamma)$ is a factorization of $K_V$ into $4$-stars.
Denote by 
$K_{U,W}$ the complete bipartite graph whose parts are $U=\mathbb{Z}\times\{0\}$ and $W=\mathbb{Z}\times \{1\}$, and consider the $1$-factor $I=\{\{i_0, i_1\}:i\in\mathbb{Z}\}$ of $K_{U,W}$.
Clearly, $K_V$ decomposed into $K_U$, $K_W  \cup I$ and $K_{U,W}\setminus I$.
One can check that
\begin{itemize}
\item $Orb_{\mathbb{Z}}(\Gamma_1)$ decomposes $K_U$,
\item $Orb_{\mathbb{Z}}(\Gamma_2)$ decomposes $K_W  \cup I$, and
\item $Orb_{\mathbb{Z}}(\Gamma_3\cup\Gamma_4)$ decomposes $K_{U,W}\setminus I$.
\end{itemize}
Hence, $\mathcal{G}$ is a decomposition of $K_V$. Considering that
the $\Gamma_i$s are pairwise vertex-disjoint and their vertex-sets partition $V$, we have that 
$\Gamma$ and each of its translates (under the action of $\mathbb{Z}$) 
are factors of $K_V$ isomorphic to a $4$-star. Therefore,
$\mathcal{G}$ is a factorization  of $K_V$ into $4$-stars.

To deal with the case $k\geq 5$, it is enough to replace the component $\Gamma_1$ of $\Gamma$ with a 
$(k-3)$-star $\Gamma'_1$ satisfying the following conditions:
\begin{align}
  V(\Gamma'_1)=V(\Gamma_1),\;\;\; \text{and} \label{replacement1}\\
  \text{$Orb_{\mathbb{Z}}(\Gamma'_1)$ decomposes $K_U$.}\label{replacement2}
\end{align}
Indeed, letting $\Gamma' = (\Gamma\setminus\Gamma_1)\ \cup\  \Gamma'_1$,  by condition \eqref{replacement1} we have that $\Gamma'$ is a $k$-star with vertex-set $V$.  
Recalling that $Orb_{\mathbb{Z}}(\Gamma_1)$ decomposes $K_U$, 
by condition \eqref{replacement2} it follows that 
$Orb_{\mathbb{Z}}(\Gamma')$ and $Orb_{\mathbb{Z}}(\Gamma)$ decompose the same graph, that is, $K_V$.
Hence, $Orb_{\mathbb{Z}}(\Gamma')$  is a factorization of $K_V$ into $k$-stars.

Let $k=h+3$ with $h\geq 2$. 
It is left to construct an $h$-star $\Gamma'_{1,h}$ satisfying conditions
\eqref{replacement1} and \eqref{replacement2}, for every $h\geq 2$.
For sake of clarity, in the rest of the proof we identify $U=\mathbb{Z}\times\{0\}$ with $\mathbb{Z}$.
Therefore, $\Gamma_1$ is the $1$-star centered in $0$ with 
$\Gamma_1(0)=\{i: i\geq 1\}$.

Let $\Delta_j$ and $\Delta^*_{j}$ be the $1$-stars centered in $c_j=2(2^j-1)$
such that
\begin{align*}
   &  \Delta_j(c_j) = \{c_j+i : 0 < i\equiv 2^j \pmod{2^{j+1}}\},\\
  & \Delta^*_j(c_j) = \{c_j+i : 0 < i\equiv 0 \pmod{2^{j}}\},
\end{align*}
for $j\geq 0$, 
and set 
$\Gamma'_{1,h} = \Delta_0\ \cup\  \Delta_1\ \cup \ldots \cup\ \Delta_{h-2}\ \cup\ \Delta^*_{h-1}$ for $h\geq 2$.
\begin{figure}
\centering
%
\begin{minipage}{15cm}
\hspace{-1.5cm}

\begin{tikzpicture}[x=2cm,y=1cm,scale=0.51]
\coordinate  (Gamma) at (-1, 5); 
\draw (Gamma)   node[above=0pt] {$\Gamma_1$};
\coordinate  (Gamma1) at (-1, 0); 
\draw (Gamma1)   node[below=0pt] {$\Gamma'_{1,2}$};
\coordinate  (Delta0) at (8, -1.5); 
\draw (Delta0)   node[below=0pt] {$\Delta_0$};
\coordinate  (Delta1star) at (9, 1.5); 
\draw (Delta1star)   node[above=0pt] {$\Delta^*_1$};
\foreach \j in {0,1,2,3,4,5,6,7,8} 
{
  \coordinate  (a\j) at (0+\j, 5); 
  \coordinate  (b\j) at (0+\j, 0);
}; 
\foreach \j in {9,10} 
{
  \coordinate  (a\j) at (0+\j, 5); 
  \coordinate  (b\j) at (0+\j, 0);
}; 
\foreach \j in  {0,1,2,3,4,5,6,7,8} 
{
  \draw (a\j) circle (4pt)[fill]  node[below=5pt] {$\j_0$}; 
};  
\foreach \j in  {0,1,2,4,6,8} 
{
  \draw (b\j) circle (4pt)[fill]  node[above=5pt] {$\j_0$}; 
};  
\foreach \j in  {3,5,7} 
{
  \draw (b\j) circle (4pt)[fill]  node[below=5pt] {$\j_0$}; 
};

\foreach \j in {1,2,3,4,5,6,7,8} 
{
    \ADD{\j}{3}{\tempB}; 
    \MULTIPLY{2}{\tempB}{\solo};     
   \MULTIPLY{-}{\solo}{\tempC}
   \ADD{\tempC}{180}{\soli};
   \draw[very thick, in=\soli, out=\solo]  (a0) to (a\j);
} 
\draw[dashed, very thick, in=156, out=24]  (a0) to (a9);
\draw[dashed, very thick, in=154, out=26]  (a0) to (a10);

\foreach \j in {4,6,8} 
{
    \ADD{\j}{3}{\tempB}; 
    \MULTIPLY{3}{\tempB}{\solo};     
   \MULTIPLY{-}{\solo}{\tempC}
   \ADD{\tempC}{180}{\soli};
   \draw[very thick, in=\soli, out=\solo]  (b2) to (b\j);
} 

\foreach \j in {1,3,5,7} 
{
    \ADD{\j}{3}{\tempB}; 
    \MULTIPLY{3}{\tempB}{\solo};     
   \MULTIPLY{-}{\solo}{\tempC}
   \ADD{\tempC}{180}{\soli};
   \draw[very thick, in=-\soli, out=-\solo]  (b0) to (b\j);
} 
\draw[dashed, very thick, in=-155, out=-35]  (b0) to (b9);
\draw[dashed, very thick, in=155, out=35]  (b2) to (b10);

%
%
\end{tikzpicture}

\end{minipage}
\label{fig2}
\end{figure}
It is not difficult to check that 
$\{\Delta_j-c_j: 0\leq j\leq h-2\}\ \cup\ \{\Delta^*_{h-1} - c_{h-1}\}$ decomposes $\Gamma_1$.
Therefore, $Orb_{\mathbb{Z}}(\Gamma'_{1,h})$ and $Orb_{\mathbb{Z}}(\Gamma_1)$ decompose the same graph, that is, $K_U$. Hence, $\Gamma'_{1,h}$ satisfies condition \eqref{replacement2}.

We show that $\Gamma'_{1,h}$ is an $h$-star satisfying condition \eqref{replacement1} by induction on $h$.
If $h=2$, then $V(\Delta_1)=\{0,1,3,5,\ldots\}$ and $V(\Delta_2)=\{2,4,6,\ldots\}$
Therefore, $\Gamma'_{1,2} = \Delta_0\ \cup\ \Delta^*_1$ 
is a $2$-star with the same vertex-set as  $\Gamma_1$.
Now assume that $\Gamma'_{1,h}$ is an $h$-star  satisfying condition \eqref{replacement1} 
for some $h\geq 2$. Recalling the definition of $\Gamma'_{1,h}$ and $\Gamma'_{1,h+1}$, 
and considering that the vertex-sets of $\Delta_{h-1}$ and $\Delta^*_{h}$ partition
$V(\Delta^*_{h-1})$, we have that $\Gamma'_{1,h+1}$ is an $(h+1)$-star with the same
vertex-set as $\Gamma'_{1,h}$, that is, $V(\Gamma_1)$, and this concludes the proof.
\end{proof}

Propositions  \ref{k1}, \ref{k2} and \ref{k>3} leave open $FP(S_k)$ only when $k=3$.
In this case, an approach similar to Theorem \ref{k>3} cannot work, as shown in the following.

\begin{Proposition} There is no $3$-star $\Gamma$ 
with vertex-set $V=\mathbb{Z}\times \{0,1, \ldots,k\}$ such that
the $\mathbb{Z}$-orbit of $\Gamma$ is an $S_3$-factorization of $K_V$
\end{Proposition}
\begin{proof} Assume for a contradiction that there exists a 
$3$-star $\Gamma$ with vertex-set $V=\mathbb{Z}\times \{0,1, \ldots,k\}$ such that
$\mathcal{G} = Orb_{\mathbb{Z}}(\Gamma)$ is a factorization of $K_V$.

We first notice that $\Gamma$ must have at least a center in $\mathbb{Z}\times \{i\}$, for every
$i\in\{0,1, \ldots,k\}$. Indeed, if $\Gamma$ has no center in $\mathbb{Z}\times \{i\}$ for some $i\in\{0,1, \ldots,k\}$, then no edge of $K_{\mathbb{Z}\times \{i\}}$ can be covered by $\mathcal{G}$. Since $\Gamma$ has $3$ centers, it follows that $k\leq 2$.
Note that if $k=2$, the centers of $\Gamma$ must be $x_0, y_1, z_2$ for some $x,y,z\in\mathbb{Z}$, but in this case the edge $\{x_0, y_1\}$ cannot lie in any translate of $\Gamma$. 
Therefore $k\leq 1$.

If $k=1$, without loss of generality we can assume that the centers of $\Gamma$
are $0_0$, $x_1$ and $y_1$ with $x\not=y$. Since the edge $\{0_0,x_1\}$ does not belong to $\Gamma$, 
it lies in some of its translates, say $\Gamma+z$ with $z\neq 0$. This is equivalent to saying that
$\{(-z)_0,(x-z)_1\}\in \Gamma$. It follows that $x-z=y$, hence $\{(y-x)_0,y_1\}\in \Gamma$.
Similarly, we can show that $\{(x-y)_0,x_1\}\in \Gamma$. It follows that $\Gamma$ cannot contain
the edges $\{0_0, (x-y)_0\}$ and $\{0_0, (y-x)_0\}$. This implies that no edge of the form
$\{w_0, (x-y+w)_0\}$ lie in any translate of $\Gamma$, contradicting again the assumption that
$\mathcal{G}$ is a factorization of $K_V$. Therefore $k=0$.

Let $V=\Z$ and denote by 
$\Delta \Gamma$ the multiset of all differences $y-x$ between any two adjcent vertices $x$ and $y$ of $\Gamma$, with $x<y$: 
\[
  \Delta \Gamma = \{y-x: \{x,y\}\in E(\Gamma),\; x<y\}.
\]
It is not difficult to see that $\mathcal{G} = Orb_{\mathbb{Z}}(\Gamma)$ is a factorization of $K_\Z$ if and only if $\Delta \Gamma = \mathbb{N}\setminus\{0\}$.
Denoting by $\Gamma+i$ the translate of $\Gamma$ obtained by replacing each vertex 
$x\in V(\Gamma)$  with $x+i$, one can easily see that $\Delta (\Gamma+i) = \Delta \Gamma$ for every $i\in\Z$.
Therefore, up to a translation, we can assume that 
the centers of $\Gamma$ are $0, x, n$ with $0<x<n$. Now, for every $i\geq n$,
denote by $\Gamma_i$ the induce subgraph of $\Gamma$ with vertex-set $\{0,1, \ldots,i\}$.
Also, let $\Gamma^*$ be the induced subgraph of $\Gamma$ on the vertices 
$\{-3,-2,-1, 0, x, n\}$. 
Clearly, 
$|\Delta \Gamma^*|=3$, $|\Delta \Gamma_i| = i-2$ and $\Delta \Gamma_i \subset \{1,2,\ldots, i\}$. Also,
since the multiset $\Delta \Gamma$ contains all positive integers with no repetition,
it follows that $\Delta \Gamma^*$ and $\Delta \Gamma_i$ are disjoint, hence
$\Delta \Gamma_i \subset \{1,2,\ldots, i\}\setminus \Delta \Gamma^*$ for every $i\geq n$. 
Then, for $i=max (\Delta \Gamma^*)$, we obtain the following contradiction:
$i-2 = |\Delta \Gamma_i| \leq |\{1,2,\ldots, i\}\setminus \Delta \Gamma^*| = i-3$.

\end{proof}
\section{The resolvability problem}
Theorem \ref{main} allows us to construct decompositions of $K_\aleph$ into $\aleph$ graphs of specified type.
More precisely, we have the following.

\begin{Corollary}\label{maincorollary}
Let $\mathcal{F}=\{F_{\alpha}: \alpha\in \mathcal{A}\}$ be an infinite family of (non-empty)
$\aleph$-bounded graphs, 
 where $\aleph=|\mathcal{A}|$. 
Then there exists a decomposition 
$\mathcal{G}=\{\Gamma_{\alpha}: \alpha\in \mathcal{A}\}$ of $K_{\aleph}$ such that
each $\Gamma_{\alpha}$ is isomorphic to $F_{\alpha}$.

Furthermore, if the domination number of some graph $F_\beta$ is less than $\aleph$, then  
$|V(K_\aleph) \setminus V(\Gamma_\beta)|= \aleph$. 
Otherwise, for every $0\preceq \aleph' \preceq \aleph$,
the decomposition $\mathcal{G}$ can be constructed so that 
$|V(K_\aleph) \setminus V(\Gamma_\beta)|= \aleph'$.
\end{Corollary}
\begin{proof} For every $\alpha\in\mathcal{A}$, set $\aleph_{\alpha}=\aleph$ if the domination number of $F_\alpha$ is less than $\aleph$; otherwise, let $0\preceq \aleph_\alpha \preceq \aleph$.
By adding to each graph $F_\alpha$ a set of $\aleph_\alpha$  isolated vertices
we obtain a graph $F'_\alpha$ whose order and domination number are $\aleph$.
Since the assumptions of Theorem \ref{main} are satisfied, there exists a factorization
$\mathcal{G'}=\{\Gamma'_{\alpha}: \alpha\in \mathcal{A}\}$ of $K_{\aleph}$ such that
each $\Gamma'_{\alpha}$ is isomorphic to $F'_{\alpha}$. By replacing $\Gamma'_\alpha$ with the isomorphic copy of $F_{\alpha}$, we obtain the desired decomposition $\mathcal{G}$.
\end{proof}

Inspired by \cite{DaHoWe14}, we ask under which conditions a decomposition $\mathcal{G}$ of $K_{\aleph}$
is resolvable, namely, its graphs can be partitioned into factors of $K_{\aleph}$, also called 
\emph{resolution classes}. It follows that a resolvable decomposition $\mathcal{G}$ of $K_{\aleph}$
must satisfy the following two conditions:
\begin{enumerate}
  \item[N1.] \label{resnec1} 
  if $\Gamma\in\mathcal{G}$ is not a factor of $K_{\aleph}$, 
  then  $|V(K_{\aleph})\setminus V(\Gamma)|\geq \min\{|\Gamma|: \Gamma\in \mathcal{G}\}$.
  \item[N2.] \label{resnec2} 
  for every $x,y,z\in V(K_\aleph)$, 
  \[\mathcal{G}(z) \subseteq \mathcal{G}(x) \cup \mathcal{G}(y)
  \;\;\Rightarrow\;\;
    \mathcal{G}(z) \supseteq \mathcal{G}(x) \cap \mathcal{G}(y),
   \]   
  where $\mathcal{G}(v)=\{\Gamma\in \mathcal{G}: v\in V(\Gamma)\}$ is the set of all graphs of 
  $\mathcal{G}$ passing through $v$.
\end{enumerate}

In the following, we easily construct decompositions of $K_\aleph$ that do not satisfy the above conditions, and therefore they are non-resolvable. 

\begin{Example} 
Let $\mathcal{F}=\{F_{\alpha}: \alpha\in \mathcal{A}\}$ be an infinite family of (non-empty)
$\aleph$-bounded graphs, 
where $\aleph=|\mathcal{A}|$. Also, assume that the domination number of at least one of its graphs,
say $F_\beta$, is $\aleph$. 
Then, by applying Corollary \ref{maincorollary} 
with $\aleph'\prec \min\{|\Gamma_\alpha|: \alpha\in \mathcal{A}\}$, we construct a decomposition that does not satisfy condition N1.

For instance, if $\aleph=|\mathbb{N}|$, each $F_{\alpha}$ is a 
countable locally finite graph
(hence, its domination number is $\aleph$) and $\aleph'=1$ for every $\beta \in \mathbb{N}$, then
we construct a decomposition 
$\mathcal{G}=\{G_\beta: \beta\in \mathbb{N}\}$ of $K_{\mathbb{N}}$ into connected regular graphs where $V(G_\beta)=\mathbb{N}\setminus \{x_\beta\}$ for some $x_\beta\in \mathbb{N}$. Clearly, no graph of 
$\mathcal{G}$ is a factor  of $K_{\mathbb{N}}$, and any two graphs of $\mathcal{G}$ have common vertices. Therefore, $\mathcal{G}$ is not resolvable.
\end{Example}

\begin{Example}
Let $\mathcal{G}$ be any decomposition of the infinite complete graph $K_{V}$ (for example, one of those constructed by Corollary \ref{maincorollary}).  
Let $y$ and $z$ be vertices not belonging to $K_{V}$ and 
set $W=V \ \cup\ \{y,z\}$. We can easily extend $\mathcal{G}$ to a non-resolvable decomposition $\mathcal{G'}$ of $K_W$ in the following way.

Choose $x\in V$ and let $\mathcal{C}$ be the following family of paths of length 1 or 2: 
\[
  \mathcal{C}=\{[y,v,z]: v\in V\setminus\{x\}\}\ \cup\ \{[x,z,y], [x,y]\}.
\]
Clearly, $\mathcal{C}$ decomposes $K_W \setminus K_V$,
hence $\mathcal{G}' = \mathcal{G}\ \cup\ \mathcal{C}$ is a decomposition of $K_{W}$. 
Also, $x,y$ and $z$ do not satisfy condition N2, since
 $\mathcal{G}'(z)\subseteq \mathcal{G}'(x)\ \cup\ \mathcal{G}'(y)$, 
while $[x,y]$ belongs to $\mathcal{G}'(x)\cap \mathcal{G}'(y)$, but not to $\mathcal{G}'(z)$.
 Therefore, $\mathcal{G}'$ is non-resolvable.  Indeed, any resolution class of $\mathcal{G}'$
 could cover the vertex  $z$ only with graphs passing through $x$ or $y$. This means that 
 the graph $[x,y]$ cannot belong to any resolution class of $\mathcal{G}'$.
\end{Example}

The following result provides sufficient conditions for a decomposition $\mathcal{G}$ to be resolvable. 
\begin{Theorem}\label{resolvableD}
Let $\mathcal{G}$ be a decomposition of the infinite complete graph $K_{\aleph}$ satisfying the following properties for some $\aleph'\prec \aleph$:
\begin{enumerate}
\item[R1.] each graph in $\mathcal{G}$ is $\aleph'$-bounded;
\item[R2.] $|\mathcal{G}(x)\ \cap \ \mathcal{G}(y)|\preceq\aleph'$ 
for every distinct $x,y\in V(K_{\aleph})$. 
\end{enumerate}
Then $\mathcal{G}$ is resolvable.
\end{Theorem}
\begin{proof} Let $\mathcal{G}=\{G_{\alpha}: \alpha\in \mathcal{A}\}$.
We consider a well-order $\prec$ on $\mathcal{A}$ 
satisfying Corollary \ref{goodgoodorder}.
Since the graphs of $\mathcal{G}$ are $\aleph'$-bounded, we have that $|\mathcal{A}|=\aleph$ and we can assume $V(K_{\aleph})= \mathcal{A}$. 
Here we need to construct an ascending chain $(\mathcal{G}_{\gamma})_{\gamma\in \mathcal{A}}$
of families 
$\mathcal{G}_{\gamma}:=\{\Gamma_{\alpha}^{\gamma}: \alpha\in \mathcal{A}_{\preceq \gamma}\}$
(where $\Gamma_{\alpha}^{\gamma'}$ is a subgraph of $\Gamma_{\alpha}^{\gamma}$ whenever 
$\gamma'\preceq \gamma$) that satisfy the following proprieties:
\begin{enumerate}
\item[($1_{\gamma}$)] each $\Gamma_{\alpha}^{\gamma}$ is a vertex-disjoint union of graphs of 
$\mathcal{G}$;

\item[($2_{\gamma}$)] for every $\alpha \in \mathcal{A}_{\preceq \gamma}$, $\gamma\in V(\Gamma_{\alpha}^{\gamma})$;

\item[($3_{\gamma}$)] $G_\gamma$ is contained in exactly one $\Gamma_{\alpha}^{\gamma}$ where $\alpha\in \mathcal{A}_{\preceq \gamma}$;

\item[($4_{\gamma}$)] for every $\alpha \in \mathcal{A}_{\preceq \gamma}$, $\Gamma_{\alpha}^{\gamma}$ is either a finite graph or  $(\aleph'\cdot|\mathcal{A}_{\preceq \gamma}|)$-bounded.
  \end{enumerate}
 
The desired resolution of $K_\aleph$ is then
$\mathcal{R}=\{\Gamma_\alpha: \alpha\in \mathcal{A}\}$, 
where $\Gamma_\alpha =\bigcup_{\gamma\in\mathcal{A}} \Gamma^\gamma_\alpha$ 
for every $\alpha\in\mathcal{A}$.
Indeed, due to properties ($1_\gamma$) and ($2_\gamma$), each $\Gamma_{\alpha}$ is a resolution class of $\mathcal{G}$ and, by property ($3_\gamma$), $\mathcal{R}$ is a partition of $\mathcal{G}$ into resolution classes.

We proceed by transfinite induction on $\gamma$.

BASE CASE. Let $0=\min  X$.
By condition $R2$, if $0$ is not a vertex of $G_0$, $|\mathcal{G}(0)\cap \mathcal{G}(x)|\preceq \aleph'$ for any $x\in V(G_0)$. Since, due to condition $R1$, $|\mathcal{G}(0)|=\aleph$, there exists $G\in \mathcal{G}(0)$ disjoint from $G_0$. 
Therefore we can define $\mathcal{G}_0=\{\Gamma_{0}^{0}\}$ where $\Gamma_{0}^{0}$ is either $G_0\cup G$ or, if $0$ belongs to $V(G_0)$, $G_0$.\\

TRANSFINITE INDUCTIVE STEP. 
For every $\gamma'\prec \gamma$, we assume there is a family 
$\mathcal{G}_{\gamma'}$ satisfying  $(i_{\gamma'})$ for $1\leq i\leq 4$. 
We show that $\mathcal{G}_{\gamma'}$ can be extended to a family 
$\mathcal{G}_{\gamma}$ that satisfies the same properties, $(i_{\gamma})$ for $1\leq i\leq 4$. 

We are going to define, recursively, the graphs $\Gamma_{\alpha}^{\gamma}$ whenever $\alpha\preceq \gamma$. 
First, we consider the case $\alpha\prec \gamma$. We start by setting $\Gamma_{\alpha}^{\prec \gamma}:=\bigcup_{\gamma'\prec\gamma} \Gamma_{\alpha}^{\gamma'}$. 
Note that  property ($4_{\gamma'}$) guarantees that $\Gamma_{\alpha}^{\prec \gamma}$ is either finite or $|\Gamma_{\alpha}^{\prec\gamma}|\preceq\aleph'\cdot |\mathcal{A}_{\preceq\gamma}|$; 
hence, $\Gamma_{\alpha}^{\prec \gamma}$ is $\aleph$-small.
\begin{itemize}
\item Base case.
If $\gamma\in V(\Gamma_{0}^{\prec\gamma})$, set $\Gamma_{0}^{\gamma}=\Gamma_{0}^{\prec\gamma}$.

If $\gamma\not \in V(\Gamma_{0}^{\prec\gamma})$, by condition R2 we have 
$|\mathcal{G}(\gamma)\cap \mathcal{G}(x)|\preceq\aleph'$ 
for every $x\in V(\Gamma_0^{\prec\gamma})$.
Since $\Gamma_0^{\prec\gamma}$ is $\aleph$-small, this means that the family of graphs of $\mathcal{G}(\gamma)$ that intersect $V(\Gamma_0^{\prec\gamma})$ is $\aleph$-small.

Moreover, any $\Gamma_{\alpha}^{\prec \gamma}$ is either finite or 
$(\aleph'\cdot |\mathcal{A}_{\preceq\gamma}|)$-bounded 
(note that $\aleph'\cdot |\mathcal{A}_{\preceq\gamma}| \prec \aleph$, 
since $|\mathcal{A}_{\preceq\gamma}|\prec \aleph$).
Hence, the set of graphs in $\mathcal{G}(\gamma)$ that are contained in some $\Gamma_{\alpha}^{\prec \gamma}$ is $\aleph$-small.

Finally, by condition R1, we have that $|\mathcal{G}(\gamma)|=\aleph$. 
Therefore, there exists a graph $G\in \mathcal{G}(\gamma)$ that is not contained in any $\Gamma_{\alpha}^{\prec \gamma}$ and such that $V(G)\cap V(\Gamma_0^{\prec\gamma})=\emptyset$.
Then, we set $\Gamma_{0}^{\gamma}=\Gamma_0^{\prec\gamma}\cup G$.

\item Inductive step. 
Let $\alpha\prec \gamma$. 
If $\gamma\in V(\Gamma_{\alpha}^{\prec\gamma})$, 
set $\Gamma_{\alpha}^{\gamma}=\Gamma_{\alpha}^{\prec\gamma}$. 
Otherwise, by proceeding as in the previous case, we obtain the existence 
of a graph $G\in \mathcal{G}(\gamma)$ that is not in any $\Gamma_{\alpha'}^{\prec \gamma}$ 
or any $\Gamma_{\alpha''}^{\gamma}$ (where $\alpha'\prec \gamma$ and $\alpha''\prec \alpha$),
and  such that $V(G)\cap V(\Gamma_{\alpha}^{\prec\gamma})=\emptyset$.
In this case, we  set $\Gamma_{\alpha}^{\gamma}=\Gamma_{\alpha}^{\prec\gamma}\cup G$.
\end{itemize}
It is left to define $\Gamma_{\gamma}^{\gamma}$. We proceed by constructing, 
recursively, an ascending chain of graphs $\Gamma_{\gamma}^{\alpha}$, for $\alpha\in \mathcal{A}_{\preceq \gamma}$, that are either finite or $(\aleph'\cdot|\mathcal{A}_{\preceq\gamma}|)$-bounded.  \begin{itemize}
\item Base case. Let us first suppose that $G_{\gamma}$ is not contained in any $\Gamma_{\alpha'}^{\gamma}$ (where $\alpha'\prec\gamma$).
Again, by conditions R1 and R2, there exists $G\in \mathcal{G}(0)$ that is also not contained in any $\Gamma_{\alpha'}^{\gamma}$ such that $G$ is either $G_{\gamma}$ or is disjoint from $G_{\gamma}$. We set $\Gamma_\gamma^0$ to be $G_{\gamma}\cup G$.
Otherwise, we set $\Gamma_\gamma^0$ to be any graph $G$ in $\mathcal{G}(0)$ that is  not contained in any $\Gamma_{\alpha'}^{\gamma}$. 
\item Inductive step. Let us suppose that $\alpha\not=0$ and that we have defined $\Gamma_{\gamma}^{\alpha'}$ for every $\alpha'\prec \alpha$. Here we set $\Gamma_{\gamma}^{\prec \alpha}$ to be $\bigcup_{\alpha'\prec \alpha}  \Gamma_{\gamma}^{\alpha'}$. Note that, for construction, $\Gamma_{\gamma}^{\prec \alpha}$ is either a finite graph or $|\Gamma_{\gamma}^{\prec \alpha}|\preceq\aleph'\cdot |\mathcal{A}_{\preceq\gamma}|$.
If $\alpha$ belongs to $V(\Gamma_{\gamma}^{\prec \alpha})$, we set $\Gamma_{\gamma}^{\alpha}$ to be $\Gamma_{\gamma}^{\prec \alpha}$.
Otherwise, proceeding as in the previous case, we obtain that there exists $G\in \mathcal{G}(\alpha)$ disjoint from $\Gamma_{\gamma}^{\prec \alpha}$ that does not belong to any of the $\Gamma_{\alpha'}^{\gamma}$. Now we set $\Gamma_{\gamma}^{\alpha}$ to be $G\cup \Gamma_{\gamma}^{\prec \alpha}$.
\end{itemize}
Then the family $\mathcal{G}_{\gamma}=\{\Gamma_{\alpha}^{\gamma}: \alpha\in \mathcal{A}_{\preceq \gamma}\}$ satisfies the properties $(1_{\gamma})$, $(2_{\gamma})$, $(3_{\gamma})$ and $(4_{\gamma})$ for construction.
\end{proof}

\begin{Remark}
A cardinal $\aleph$ is said to be \textit{regular} if any $\aleph$-small union of $\aleph$-small sets (resp. graphs) is still an $\aleph$-small set (resp. graph) otherwise it is said to be \textit{singular}.
It is easy to see that, for regular cardinals, conditions $R1$ and $R2$ of Theorem \ref{resolvableD} can be relaxed to:
\begin{enumerate}
\item[$R1'.$] each graph in $\mathcal{G}$ is $\aleph$-small;
\item[$R2'.$] $|\mathcal{G}(x)\ \cap \ \mathcal{G}(y)|\prec\aleph$ 
for every distinct $x,y\in V(K_{\aleph})$. 
\end{enumerate}
However, if $\aleph$ is a singular cardinal, then conditions $R1'$ and $R2'$ are no longer sufficient.
Indeed, we can construct a decomposition $\mathcal{G}$ of $K_{\aleph}$ into $\aleph$-small graphs such that
\begin{enumerate} 
  \item[a.] $|\mathcal{G}|$ is $\aleph$-small,
  \item[b.] $\mathcal{G}$ satisfies conditions $R1'$ and $R2'$,
  \item[c.] there are two (possibly isolated) vertices $x$ and $y$ belonging to every graphs of $\mathcal{G}$, that is, $\mathcal{G}=\mathcal{G}(x) \cap \mathcal{G}(y)$
\end{enumerate}
 Then, choosing any vertex $z$ such that $\mathcal{G}(z)\not=\mathcal{G}$, we have that
  \[\mathcal{G}(z) \subseteq \mathcal{G}(x) \cup \mathcal{G}(y)=\mathcal{G}
  \;\;\ \mbox{ but }\;\;
    \mathcal{G}(z) \not\supseteq \mathcal{G}(x) \cap \mathcal{G}(y)=\mathcal{G}.
   \]   
This means that condition N2 does not hold, therefore the decomposition $\mathcal{G}$ is not resolvable.
\end{Remark}
%

We conclude by showing that there is always a resolution for an `almost' $2$-design with
blocks that are $\aleph'$-bounded for some $\aleph'\prec \aleph$,  
that is, a decomposition of $K_{\aleph}$ whose graphs are almost all  $\aleph'$-bounded complete graphs. 
This extends some results on the resolvability of $2$-designs given in \cite{DaHoWe14}.

\begin{Proposition}
  Let $\mathcal{G}$ be a decomposition of the infinite complete graph $K_\aleph$ into $\aleph'$-bounded graphs for some $\aleph'\prec \aleph$, where $\aleph'$ is not necessarily infinite. 
  If the subset of $\mathcal{G}$ consisting of all non-complete graphs is $\aleph'$-bounded,
  then $\mathcal{G}$ has a resolution.
\end{Proposition}
\begin{proof}
 By assumption, condition R1 of Theorem \ref{resolvableD} holds. To prove that $\mathcal{G}$ satisfies condition R2 for some $\aleph''\prec \aleph$, we assume for a contradiction the existence of vertices $x$ and $y$ such that $|\mathcal{G}(x)\ \cap \ \mathcal{G}(y)|\succ \aleph'':=(\aleph'+1)$. It follows that there are at least two complete graphs in $\mathcal{G}(x)\ \cap \ \mathcal{G}(y)$, meaning that the edge $\{x,y\}$ is covered more than once by graphs in $\mathcal{G}$, and this is a contradiction. The assertion follows from Theorem \ref{resolvableD}.
 \end{proof}


\begin{thebibliography}{99}
\bibitem{BeuCa94}
A. Beutelspacher, P.J. Cameron,
Transfinite methods in geometry,
Bull. Belg. Math. Soc. 3 (1994), 337--347.


\bibitem{BoMa10}
S. Bonvicini and G. Mazzuoccolo,
\PaperTitle{Abelian $1$-factorizations in infinite graphs},
\JournalName{European J. Combin.}, 31: 1847--1852, 2010.


\bibitem{Ca84}
P. J. Cameron,
Infinite versions of some topics infinite geometry, in: F.C. Holroyd, R.J. Wilson(Eds.), Geometrical Combinatorics, in: Res. Notes Math. Ser., vol. 114, Pitman, Boston, 1984, pp.13--20.

\bibitem{Ca95}
P. J. Cameron, Note on large sets of infinite Steiner systems,\ 
J. Combin. Des. 3 (1995), 307--311.

\bibitem{Ca97}
P. J. Cameron, The random graph, The mathematics of Paul Erdős, II, Algorithms Combin., 14, Berlin: Springer, pp. 333-351.

\bibitem{CaWe12}
P. J. Cameron and B. S. Webb,
Perfect countably infinite Steiner triple systems,
Australas. J. Combin. 54 (2012), 273--278.

\bibitem{ChiGrGrWe10}
K. M. Chicot, T. S. Griggs, M. J. Grannell and B. S. Webb,
On sparse countable infinite Steiner triple systems,
J. Combin. Des. 18 (2010), 115--122.

\bibitem{Costa20}
S. Costa,
\PaperTitle{A complete solution to the infinite Oberwolfach problem},
J. Combin. Des. 28 (2020), 366--383.

\bibitem{DaHoWe14}
P.~Danziger, D. Horsley and B. S. Webb,
\PaperTitle{Resolvability of infinite designs},
J. Combin. Theory A 123 (2014), 73--85. 

\bibitem{Die17}
R. Diestel, Graph Theory, 5th ed., Springer-Verlag, Heidelberg, 2017.


\bibitem{Fr94}
F. Franek,
Isomorphisms of infinite Steiner triple systems,
Ars Combin. 38 (1994), 7--25.

\bibitem{GrGrPh87}
M. J. Grannell, T. S. Griggs and J. S. Phelan,
Countably infinite Steiner triple systems,
Ars Combin. 24B (1987), 189--216.

\bibitem{GrGrPh91}
M. J. Grannell, T. S. Griggs and J. S. Phelan,
On infinite Steiner systems, 
Discrete Math. 97 (1991), 199--202.

\bibitem{Ko77}
E. K\"{o}hler,
Unendliche gefaserte Steiner systeme,
J. Geom. 9 (1977), 73--77.

\bibitem{Qua80}
R. W. Quackenbush,
Algebraic speculations about Steiner systems,
Ann. Discrete Math. 7 (1980), 25-30.

\bibitem{Rado}
R. Rado, 
\PaperTitle{Universal graphs and universal functions}, 
\JournalName{Acta Arith.}, 9: 331-340, 1964.



\end{thebibliography}
\end{document}